\documentclass[a4paper,12pt]{article}
\usepackage[LCY,T2A]{fontenc}
\usepackage[cp1251]{inputenc}
\usepackage{amsthm}
\usepackage{amsmath}
\usepackage{amssymb}
\usepackage{mathrsfs}
\usepackage[english,russian]{babel}
\usepackage[active]{srcltx}

\begin{document}

\centerline {\bf Comonotone Second Jackson's Inequality}

\begin{center}
Pleshakov Mikhail Gennad'evich, Saratov State University,
Russia,\\
410012, Saratov, Astrakhanskaya st., 83, pleshakovmg@gmail.com
\end{center}

{\bf Abstract.}
Let $2s$ points $y_i=-\pi\le y_{2s}<\ldots<y_1<\pi$ be given.
Using these points, we define the points $y_i$ for all integer
indices $i$ by the equality $y_i=y_{i+2s}+2\pi$. We shall write
$f\in\bigtriangleup^{(1)}(Y)$ if $f$ is a $2\pi$-periodic function
and $f$ does not decrease on $[y_i, y_{i-1}]$ if $i$ is odd; and
$f$ does not increase on $[y_i, y_{i-1}]$ if $i$ is even. We
denote $E_n^{(1)}(f;Y)$ the value of the best uniform comonotone
approximation.  In this article the following Theorem -- the
comonotone analogue of  second Jackson's Inequality -- is
proved.\newline \textbf{Theorem.} {\it If
$f\in\bigtriangleup^{(1)}(Y)\bigcap {\Bbb W}^r$, $r>2$, then
$$E_n^{(1)}(f;Y)\le \frac c{n^r}, $$
where $c=c(r,Y)=const$ depending only on $r$ and $Y$, ${\Bbb W}^r$
Sobolev space.}


При каждом $r\in\Bbb N$ пусть $\Bbb W^r$ --  класс
$2\pi$-периодических функций, имеющих $r-1$ абсолютно непрерывную
производную и
$$|f^{(r)}(x)|\le 1,\ \text{п.в.\ на}\ \Bbb R.$$
Будем рассматривать $2\pi$-периодические функции $f\in\Bbb W^r, r\geq 2$.

Обозначим ${\Bbb T}_n, n\in\Bbb N$ - пространство тригонометрических
полиномов вида
$$ \tau_n(x)= a_0 + \sum\limits_{k=1}^n a_k\cos kx + b_k\sin kx,
\ \ x\in\Bbb R $$
порядка $\leq n$.

Пусть на промежутке $[-\pi ,\pi )$ заданы $2s$ точек $y_i$
$$-\pi\leq y_{2s}<\ldots <y_1<\pi.$$
Отправляясь от этих точек, при помощи равенства
$$y_i:=y_{i+2s}+2\pi$$
определим точки $y_i$ для всех целых индексов $i$.

Обозначим $Y:=\{ y_i\}_{i=1}^{2s}$,
 $$\Pi (x):=\prod_{i=1}^{2s} \sin\frac 12(x-y_i),$$
 и заметим, что $\Pi\in {\Bbb T}_s$, т.е. $\Pi$ - тригонометрический
полином порядка $s$.

Пусть $f$ -- непрерывно дифференцируемая функция.
Будем писать $f\in\bigtriangleup^{(1)}(Y)$, если
$$f'(x)\Pi (x)\geq 0,\ \ x\in\Bbb R.                 \eqno (1)   $$

\noindent
Положим
$$E_n^{(1)}(f;Y):=\inf\limits_{\tau\in\bigtriangleup^{(1)}(Y)\cap\Bbb T_n}
\| f-\tau \|.$$

\medskip
\noindent
В главе 2 доказывается

\noindent {\bf Т е о р е м а 1.} {\it Если
$f\in\bigtriangleup^{(1)}(Y)\cap\Bbb W^r$,
 $r\geq 2$,
то
$$E_n^{(1)}(f;Y)\leq\frac c{n^r}.                    \eqno (2)   $$
где $c=c(r,Y)$ - постоянная, не зависит от $n$ и $f$.}

\medskip
{\it З а м е ч а н и е  1.} Теорема 1 верна и в случае $r=1$, это
следует из теоремы 1 [1]. Однако функция $f\in\Bbb W^1$ не обязана
иметь непрерывную производную и поэтому в случае $r=1$ множество
$\bigtriangleup^{(1)}(Y)$ нуждается в отдельном определении.
Теорема 1 является аналогом соответствующих результатов о
комонотонном приближении непрерывных функций алгебраическими
многочленами, полученных в работах [2-20].

\medskip
Доказательство теоремы разобьем на несколько пунктов. Всюду далее через
$c_1$, $c_2$, $c_3$, ... обозначены положительные числа (константы),
которые  зависят только от $s$ или $r$ или $s$ и $r$, через $C_i=C_{i,l}$
обозначены положительные числа (константы), которые зависят только от $s$, $r$
и $l$.

\medskip
\noindent $1^o$. Некоторые обозначения.

Зафиксируем натуральное $n$ и разобьем прямую $\Bbb R$ равноотстоящими точками
$$ x_j:= -\frac {j\pi}{n},\ \ j\in\Bbb Z.$$
Обозначим
$$ I_j:=[x_j;x_{j-1}],\ \ h:=\frac {\pi}{n}.$$

\medskip
\noindent {\bf Л е м м а 1.} {\it Пусть $g\in\Bbb W^{r-1}$ и
пусть  при некотором фиксированном $j\in\Bbb Z$ на каждом из
$2r-3$ отрезков $I_j,\ldots , I_{j+2(r-2)}$ существует точка
$t_{\nu}\in I_{\nu},\ \ \nu =\overline {j,j+2(r-2)}$, такая, что
$$| g(t_{\nu}) |\leq h^{r-1},                                      \eqno (3) $$
тогда для любого $x\in\bigcup\limits_{\nu =j}^{j+2(r-2)} I_{\nu}$
справедливо неравенство
$$| g(x) |\leq c_1h^{r-1}.                                         \eqno (4) $$}

{\it Д о к а з а т е л ь с т в о}. Рассмотрим разделенную разность
порядка $r-1$ функции $g$ в узлах $x$, $t_j$, $t_{j+2}$, ...
,$t_{j+2(r-2)}$
\begin{gather*} [t_j,\ldots
,t_{j+2(r-2)},x;g]=\frac {g(t_j)}{(t_j-t_{j+2})\ldots
(t_j-t_{j+2(r-2)})(t_j-x)}+
\ldots \\
+\frac { g(t_{j+2(r-2)}) }{(t_{j+2(r-2)}-t_j)\ldots (t_{j+2(r-2)}-t_{j+2(r-3)})
(t_{j+2(r-2)}-x)} +\\
+ \frac {g(x)}{(x-t_j)\ldots (x-t_{j+2(r-2)})}.
\end{gather*}
Из условия $g\in\Bbb W^{r-1}$ леммы  и свойств разделенной
разности ( см., например, [21, с. 13]) следует оценка
$$| [t_j,\ldots ,t_{j+2(r-2)},x;g] |\leq\frac 1{(r-1)!}.           \eqno (5) $$
После тождественных преобразований получаем

\begin{gather*}
|g(x)|\leq\frac {| (x-t_j)\ldots (x-t_{j+2(r-2)} )|}{(r-1)!} + \\
+\left\vert\frac {g(t_j)(x-t_{j+1})\ldots (x-t_{j+2(r-2)})}{(t_j-t_{j+2})
\ldots (t_j-t_{j+2(r-2)}) } \right\vert+\ldots  \\
 \ldots+\left\vert\frac{g(t_{j+2(r-2)})(x-t_j) \ldots (x-t_{j+2(r-3)})}
{(t_{j+2(r-2)}-t_j)\ldots (t_{j+2(r-2)}-t_{j+2(r-3)})}\right\vert\le \\
\le h^{r-1}\left(\frac {(2r-3)^{r-1}}{(r-1)!}+(r-1)(2r-3)^{r-2}\right).
\end{gather*}
Лемма доказана.

\medskip

\noindent $2^o$. Разбиение индексов $j\in\Bbb Z$ на три типа.

Каждый из индексов $j\in\Bbb Z$ отнесем к одному и только одному из трех типов.

\noindent {\it О п р е д е л е н и е 1.} Индекс $j$ назовем
индексом первого типа и будем писать $j\in V_1$, если для любого
$x\in I_j$ справедливо неравенство
$$| f'(x)|\leq c_1h^{r-1},                                      \eqno (6) $$
обозначим
$$E_1:=\bigcup\limits_{j\in V_1}I_j.                            \eqno (7) $$

\noindent
Индекс $j$ назовем индексом второго типа
и будем писать $j\in V_2$, если он не является индексом первого типа
 и для любого $x\in I_j$ справедливо неравенство
$$| f'(x)|\geq h^{r-1}.                                    \eqno (8) $$
Обозначим
$$E_2:=\bigcup\limits_{j\in V_2}I_j.                       \eqno (9) $$

\noindent
Все остальные индексы отнесем к третьему типу и множество всех
индексов третьего типа обозначим через $V_3$.
Обозначим
$$E_3:=\bigcup\limits_{j\in V_3}I_j.                       \eqno (10) $$

\medskip

\noindent {\bf Л е м м а 2.} {\it Индексов третьего типа не может
быть более $2r-4$ штук подряд.}

{\it Д о к а з а т е л ь с т в о}. Допустим от противного, что
имеется $2r-3$ подряд индексов третьего типа $j, j+1,\ldots
,j+2(r-2)$. Это означает, что на каждом отрезке $I_{\nu},\ \ \nu
=\overline {j, j+2(r-2)}$ найдется по крайней мере одна точка
$t_{\nu}$, для которой справедливо неравенство $$|
f'(t_{\nu})|\leq h^{r-1}.$$ По лемме 1 это значит, что для любого
$x\in\bigcup\limits_{\nu =j}^{j+2(r-2)} I_{\nu}$ справедлива
оценка
 $$| f'(x)|\leq c_1h^{r-1}, $$ а это противоречие, т.к. на множествах
третьего типа должны быть точки $y$, в которых $$| f'(y)|>c_1h^{r-1}.$$
Лемма доказана.  \medskip

\noindent $3^o$. Построение множества $\Omega$.

Каждый набор индексов $j+1,\ldots ,j+\mu$, состоящий из не менее семи (т.е.
$\mu\geq  7$) следующих подряд индексов второго типа  назовем пачкой, если
$j\overline {\in} V_2$ и $(j+\mu +1)\overline {\in} V_2$.

Если не существует ни одной пачки, то положим $W_1:=\Bbb Z$. Если пачки
существуют, то каждый индекс $j$, не входящий ни в одну из пачек
отнесем к одной из следующих двух групп : $W_1$ или $W_2$.

Если между двумя соседними пачками не содержится ни одного индекса
первого типа, то все индексы между этими двумя пачками отнесем к группе
$W_2$; в противном случае все индексы между этими двумя пачками отнесем
к группе $W_1$.

Положим
$$ M:=\bigcup\limits_{j\in W_1} I_j.$$

Для любого множества $X\subset\Bbb R$ обозначим
$$X^*:=\bigcup\limits_{j: I_j\cap X\ne 0} I_j,\ \ X^{**}:=\left( X^*\right)^*,
\ \ X^{***}:=\left( X^{**}\right)^*,$$
и т.д.
Положим $$\Omega :=M^{***}.$$

\medskip
\noindent $4^o$. Оценка производной на $\Omega $.

\medskip
\noindent {\bf Л е м м а 3.} {\it Для всех $x\in\Omega $ имеет
место оценка
$$| f'(x) |\leq c_2h^{r-1}.                           \eqno (11) $$}

{\it Д о к а з а т е л ь с т в о}. Зафиксируем $j\in\Bbb Z$ такое, что
$x\in I_j$ и рассмотрим два множества
$$\overline {A}_j :=\bigcup\limits_{\nu =j}^{j+14r} I_{\nu} \ \ \
и\ \ \ \underline {A}_j :=\bigcup\limits_{\nu =j-14r}^{j} I_{\nu}.$$

\noindent
1. Пусть множество $\overline {A}_j$ пересекается с некоторой пачкой
и множество $\underline {A}_j$ также пересекается с некоторой пачкой. Тогда
по построению множества $\Omega$ среди индексов
$\nu =\overline {j-14r, j+14r}$ найдется по крайней мере один индекс первого
типа, который обозначим через $\nu ^*$.
Согласно определению индексов первого типа для любого $y\in I_{\nu ^*}$
имеем
$$| f'(y) |\leq c_1h^{r-1}.$$
Теперь рассмотрим разделенную разность порядка $r-1$

\noindent
\ $[t_{j_0},\ldots ,t_{j_{r-2}},x;f']$
в узлах $x$ и $t_{j_k}\in I_{\nu^*}$, где $t_{j_k}:=x_{\nu^*}+\frac {\pi k}
{n(r-2)}$,\ $k=\overline {0, r-2}$.

Повторяя рассуждения леммы 1, получим оценку
$$| f'(x) |\leq c_2h^{r-1},$$
т.е. в случае $1$ оценка $(11)$ доказана.

\medskip
\noindent
2. Пусть по крайней мере одно из множеств $\underline {A}_j$,
$\overline {A}_j$ не пересекается ни с одной пачкой, скажем
$\overline {A}_j$. Если по крайней мере один из индексов $j,\ldots ,j+14r$,
составляющих множество $\overline {A}_j$, является индексом первого типа,
то рассуждаем аналогично пункту 1.
В противном случае рассуждаем следующим образом.
По построению множества $\Omega$ среди индексов $j,\ldots ,j+6$ имеется по
крайней мере один индекс третьего типа; среди индексов $j+7,\ldots ,j+13$
также имеется по крайней мере один индекс третьего типа и т.д.
Другими словами, среди всех индексов $\nu =j,\ldots ,j+14r$ найдутся по
крайней мере $r-1$ штук индексов третьего типа $\nu_1,\ldots ,\nu_{r-1}$
таких, что
$$8\leq\vert \nu_i-\nu_q\vert\leq 14r.$$
Поэтому на множестве $\overline {A}_j$ существуют $r-1$ точек $t_i$,
$i=\overline {1, r-1}$ таких, что
$$| f'(t_i)|\leq h^{r-1}$$ и
$$6h\leq | t_i-t_q |\leq 14rh.$$
Рассуждая как в лемме 1, получим оценку
$$| f'(x) |\leq c_2h^{r-1}.$$
Лемма доказана.

\medskip
\noindent $5^o$. Определение ``маленькой'' функции $g_1(x)$.

При всяком $j\in\Bbb Z$ обозначим
$$S_j(x):=\frac {1}{b}\int\limits_{x_j}^{x}(t-x_j)^{r-2}(x_{j-1}-t)^{r-2}
dt,$$
$$b:=\int\limits_{x_j}^{x_{j-1}}(t-x_j)^{r-2}(x_{j-1}-t)^{r-2}dt,$$
и заметим, что
$$ S_j(x_j)=0,\ \ S_j(x_{j-1})=1;                              \eqno (12) $$
$$0\leq S_j(x)\leq 1,\ \ x\in I_j;                             \eqno (13) $$
$$S_j^{(q)}(x_j)=S_j^{(q)}(x_{j-1})=0,\ \ q=\overline {1,r-2}; \eqno (14) $$
$$\left\vert S_j^{(q)}(x)\right \vert \leq c_3h^{-q},\ \ q=\overline {1,r-1},
x\in I_j. \eqno (15) $$ Положим
\[g_1(x):=\left\{\begin{array}{rl}
f'(x),&\ \mbox {если}\ \ x\in M^*,\\
                 0    ,&\ \mbox {если}\ \ x\in\Bbb R \setminus M^{**}.
          \end{array}\right.\]                    (15')
На $M^{**}\setminus M^*$ определим $g_1(x)$ следующим образом.
Если $I_j\subset\overline {M^{**}\setminus M^*}$ и
$I_{j+2}\overline\subset \Omega$, то для $x\in I_j$ положим
$g_1(x):=f'(x)S_j(x),\ \ x\in I_j$; в противном случае (т.е.
$I_{j+2}\subset \Omega$) для $x\in I_j$ положим
$g_1(x):=f'(x)(1-S_j(x))$.

\medskip
\noindent {\bf Л е м м а 4.} {\it Функция $g_1(x)$ имеет свойства

\noindent
$$\ \ g_1\in c_4\Bbb W^{r-1}.                         \eqno (16) $$

\noindent
$$\ \ \Vert g_1 \Vert_{\Bbb R}\leq c_2h^{r-1}.        \eqno (17) $$

\noindent
$$\ \ g_1(x)\Pi (x)\geq 0,\ \ x\in\Bbb R.             \eqno (18) $$

\noindent
$$ \ (f'(x)-g_1(x))\Pi (x)\geq 0,\ \ x\in\Bbb R.      \eqno (19) $$}

{\it Д о к а з а т е л ь с т в о}. По условию теоремы $\vert
f^{(r)}(x)\vert\leq 1 $ почти всюду на $\Bbb R$; по лемме 3
$$\vert f'(x)\vert\leq c_2h^{r-1},\ \ x\in\Omega .         \eqno (20) $$
В силу $(12)$ и $(14)$ функция $g_1(x)$ имеет непрерывную
$(r-2)$-ю производную, нетрудно заметить, что эта $(r-2)$-я
производная также абсолютно непрерывна. Докажем оценку
$$| g_1^{(r-1)}(x)|\leq c_5                                    \eqno (21) $$
для почти всех $x\in\Bbb R$. Если $x\overline\in M^{**}\setminus
M^*$, то оценка $(21)$ очевидна. Пусть $x\in I_j\subset\overline
{M^{**}\setminus M^*}$. Для доказательства $(21)$ используем
неравенство типа Колмогорова
$$\Vert g^{(j)}\Vert_{[a,b]}\leq с_6\left( (b-a)^{r-j}\Vert g^{(r)}
\Vert_{[a,b]} + (b-a)^{-j}\Vert g\Vert_{[a,b]}\right),\ \ j=\overline {0,r}.$$
которое для функции $g(x)=f'(x)$ и $[a, b]=I_j$ принимает вид
$$\Vert f^{(j+1)}\Vert_{I_j}\leq c_6h^{r-j-1}\Vert f^{(r)}\Vert_{I_j} +
c_6h^{-j}\Vert f'\Vert_{I_j}\leq c_7h^{r-j-1}. \eqno (22) $$
Рассмотрим
\begin{gather*} (f'(x)S_j(x))^{(r-1)}= f^{(r)}(x)S_j(x)
+ \binom {r-1}1f^{(r-1)}(x)S_j'(x)+
\ldots \\
\ldots + f'(x)S_j^{(r-1)}(x).  \end{gather*}

Пользуясь (15) и (22) получаем
$$\left\vert\left( f'(x)S_j(x)\right)^{(r-1)}\right\vert\leq c_82^{r-1}.
                                                                \eqno (23) $$
Т.е. оценка $(21)$ доказана, а значит доказано соотношение $(16)$.
Оценка $(17)$ следует из (11) и очевидной оценки
$$\vert g_1(x)\vert\leq\vert f'(x)\vert,\ \ x\in\Omega .    $$
Из этой же оценки следует неравенство $(19)$.

Hеравенство $(18)$ очевидно, когда $x\overline\in\ M^{**}\setminus
M^*$, если же

\noindent $x\in\overline {M^{**}\setminus M^*}$, то $(18)$ следует
из $(13)$.

Лемма доказана.

\medskip
Обозначим
$$ G_1(x):=f(0) +\int\limits_0^x g_1(u)du $$
и заметим, что
$$ G_1(x)=Bx + \tilde {G_1}(x),                                \eqno (24) $$
где $B=const$, $\tilde {G_1}-2\pi$ - периодическая функция.

\medskip

\noindent $6^o$. Приближение ``ступеньки''.

Обозначим
$$J_{l,n}(t):=\frac {1}{\gamma_{l,n}}\left(\frac {\sin\frac {nt}{2}}
{\sin\frac {t}{2}}\right)^{2l},\ l\in\Bbb N, $$
-- ядро типа Джексона, в котором
$$\gamma_{l,n}:=\int\limits_{-\pi}^{\pi}\left(\frac {\sin\frac {nt}{2}}
{\sin\frac {t}{2}}\right)^{2l}dt.$$ Ядра типа Джексона обладают
следующими свойствами ([22], см. также например [23, с. 130]):

\noindent
$$\ \int\limits_{-\pi}^{\pi}J_{l,n}(t)dt=1,                     \eqno (25) $$

\noindent
$$\ J_{l,n}\in{\Bbb T}_{l(n-1)},                                \eqno (26) $$

\noindent
$$\ \frac {1}{C_9}n^{2l-1}<\gamma_{l,n}<C_{10}n^{2l-1};         \eqno (27) $$

\noindent
при любом $\delta\geq 0$ и $\nu =\overline {0,2l-2}$
$$\int\limits_{\delta}^{\pi}(1+nt)^{\nu}J_{l,n}(t)dt
\leq\frac {C_{11}}{(1+n\delta )^{2l-\nu -1}}, \eqno (28) $$ в
частности
$$\int\limits_{-\pi}^{\pi}(1+n|t|)^{\nu}J_{l,n}(t)dt\leq C_{12}.  \eqno (29) $$

\noindent {\it О п р е д е л е н и е 2.} Для каждой точки
$x^*\overline\in Y$ обозначим
$$ {T}_{l,n}(x;x^*;s;Y):=\frac {1}{d_{l,n}(x^*;s;Y)}
\int\limits_{x^* -\pi}^x J_{l,n}(t-x^*) \frac {\Pi (t)}{\Pi
(x^*)}dt,                                   \eqno (30) $$ где
$$d_{l,n}(x^*;s;Y):=\int\limits_{x^*-\pi}^{x^*+\pi}J_{l,n}(t-x^*)
\frac {\Pi (t)}{\Pi (x^*)}dt. $$

Положим
\[\chi (x;x^*):=\left\{\begin{array}{rl} 0,\ \mbox {если} \ \ x\le x^*,
\\
                       1,\ \mbox {если} \ \ x> x^*;
\end{array}\right.\]                                                       (31)
$$\delta_n(x;x^*):=\min\left\{ 1,\frac {1}{n|\sin\frac {x-x^*}{2}|}
\right\}.$$

\noindent {\bf Л е м м а 5.} {\it Если $l>s$ и
$$\min\limits_{i\in\Bbb Z}| x^* - y_i |\geq 2sC_{12}\frac {\pi}{n},
                                                              \eqno (32)  $$
то для функции $T (x)=T_{l,n}(x;x^*;s;Y)$ имеют место соотношения
$$\frac {1}{2} < d_{l,n}(x^*;s;Y) < \frac {3}{2},             \eqno (33)  $$
$$T (x)=\frac {1}{2\pi}x + \tilde {T}(x),                    \eqno (34)  $$
 где       $$\tilde T\in\Bbb T_{l(n-1)+s},                    \eqno (34') $$
$$\Pi (x^*)T '(x)\Pi (x)\geq 0,\ \ x\in\Bbb R,                \eqno (35)  $$
$$|T '(x) |\leq C_{13}n\delta_n^{2(l-s)}(x;x^*),\ \ x\in\Bbb R,
                                                              \eqno (36) $$
$$|T '(x) |\geq\frac {1}{C_{14}}n\delta_n^{2l}(x;x^*)\left\vert\frac {\Pi (x)}
{\Pi (x^*)}\right\vert , \eqno (37) $$
$$x\in\bigcup\limits_{\nu\in\Bbb Z}[x^* + \frac {2\pi\nu}{n}
+\frac {\pi}{2n}, x^* + \frac {2\pi\nu}{n} + \frac {2\pi}{n}-\frac {\pi}{2n}
];$$
$$|\chi (x;x^*)-T (x)|\leq C_{15}\delta^{2(l-s)-1}_n(x,x^*),
\ \ |x - x^*|\leq 2\pi . \eqno (38) $$}

\noindent
{\it Д о к а з а т е л ь с т в о}. Представим
$d=d_{l,n}(x^*;s;Y)$ в виде
$$d=\int\limits_{-\pi}^{\pi}J_{l,n}(t)\frac {\Pi (x^* +t)}{\Pi (x^*)} dt.$$
Заметим, что
\begin{gather*} |\Pi '(\theta )| =\frac
12\left\vert\sum\limits_{i=1}^{2s}\cos \frac {\theta -
y_i}{2}\prod\limits_{\nu =1,\nu\ne i}^{2s}\sin\frac
{\theta -y_{\nu}}{2}\right\vert =\\
=\frac 12\left\vert\sum\limits_{i=1}^{2s}\cos\frac {\theta
-y_i}{2} \prod\limits_{\nu =1, \nu\ne i}^{2s}\left( \sin\frac
{\theta -x^*}{2} \cos \frac{x^*-y_{\nu}}{2} + \cos\frac {\theta
-x^*}{2}\sin\frac
{x^*-y_{\nu}}{2}\right)\right\vert \\
\leq\frac {1}{2}\sum\limits_{i=1}^{2s}\prod\limits_{\nu =1, \nu\ne
i}^{2s} \left(\frac {|\theta -x^*|}{2}+\left\vert\sin\frac {x^* -
y_{\nu}}{2}\right \vert\right) ,\ \ \theta\in\Bbb R.
\end{gather*}

Поэтому для некоторого $\theta$, лежащего между точками $x^*$ и
$x^*+t$ , имеем
\begin{gather*} \left\vert 1 -\frac {\Pi (x^*
+t)}{\Pi (x^*)}\right\vert =\left\vert
\frac {\Pi (x^* +t)-\Pi (x^*)}{\Pi (x^*)}\right\vert =\\
=\left\vert\frac {t\Pi '(\theta )}{\Pi (x^*)}\right\vert\leq\frac
{|t|}{2} \sum\limits_{i=1}^{2s}\frac {1}{|\sin\frac {x^*-y_i}{2}|}
\prod\limits_{\nu =1,\nu\ne i}^{2s}\left(\frac {|\theta - x^*|}
{2|\sin (\frac {x^* -y_{\nu}}{2})|} + 1 \right), \end{gather*}
поскольку $|\theta - x^*|<t$, а в силу $(32)$
$$\left\vert\sin\frac {x^* -y_i}{2}\right\vert >\frac {2sC_{12}}{n}=:
\frac {C_{16}}{n},$$
то
$$\left\vert 1 -\frac {\Pi (x^* +t)}{\Pi (x^*)}\right\vert\leq
| t|\frac {n}{C_{16}}s\left(1+\frac {n| t|}{2C_{16}}\right)^{2s-1}<$$
$$<\frac {s}{C_{16}}(1+n| t| )^{2s}.$$
Следовательно, с учетом $(25)$ и $(29)$ получаем

\begin{gather*} |d-1| =\left\vert\int\limits_{-\pi}^{\pi}J_{l,n}(t)\left(
\frac {\Pi (x^* +t)}{\Pi (x^*)} -1\right)\,dt\right\vert <\\
<\frac {s}{C_{16}}\int\limits_{-\pi}^{\pi}J_{l,n}(t)(1+n| t|
)^{2s}dt \leq\frac {s}{C_{16}}C_{12}=1/2. \end{gather*}

Таким образом оценка $(33)$ доказана. Соотношения $(34)$ и $(35)$
очевидны.

Докажем оценку $(36)$. Очевидно
$$T'(x)=\frac {1}{d}J_{l,n}(x-x^*)\frac {\Pi (x)}{\Pi (x^*)}.    \eqno (39) $$
В силу определения ядра Джексона имеем

\begin{gather*}
J_{l,n}(t)\leq\frac {1}{\gamma_{l,n}}\left(\min\left\{ n,\frac
{1}{|\sin\frac {t}{2}|}\right\}\right)^{2l}=\\
=\frac {n^{2l}}{\gamma_{l,n}}\left(\min\left\{ 1,\frac {1}{n|\sin\frac
 {t}{2}|}\right\}\right)^{2l},\ \ t\in\Bbb R,\end{gather*}
т.е. согласно $(27)$
$$J_{l,n}(x-x^*)\leq C_{10}n\delta_n^{2l}(x;x^*).             \eqno (40) $$
Теперь, с учетом $(32)$ и равенства $2sC_{12}=C_{16}$,
\begin{gather*} \left\vert\frac {\sin\frac {x-y_i}{2}}{\sin\frac
{x^*-y_i}{2}}\right\vert =\left\vert\frac {\sin\frac
{x-x^*}{2}\cos\frac {x^*-y_i}{2} + \cos\frac
{x-x^*}{2}\sin\frac {x^*-y_i}{2}}{\sin\frac {x^*-y_i}{2}}\right\vert\leq\\
\leq 1+\left\vert\frac {\sin\frac {x-x^*}{2}}{\sin\frac
{x^*-y_i}{2}} \right\vert\leq 1+\frac {n|\sin\frac
{x-x^*}{2}|}{C_{16}}\leq
 2\delta_n^{-1}(x,x^*),\end{gather*}
а значит
$$|T'(x)|\leq 2^{2s+1}C_{10}n\delta_n^{2l-2s}(x;x^*)=:C_{13}n\delta_n
^{2(l-s)}(x;x^*).$$ Оценка $(36)$ доказана.

Докажем оценку $(37)$. При $$x\in\bigcup\limits_{\nu\in\Bbb
Z}\left[ x^* +\frac {2\pi\nu}{n} +\frac {\pi}{2n},x^* +\frac
{2\pi\nu}{n} +\frac {2\pi}{n} -\frac {\pi}{2n}\right]$$ имеем
$$\left\vert\sin\frac {n(x-x^*)}{2}\right\vert\geq\frac {\sqrt 2}2,$$
откуда
$$|T'(x)|\geq\frac {2^{-l}}{d\gamma_{l,n}}\left(\frac {1}{|\sin\frac {x-x^*}
{2}|}\right)^{2l}\left\vert\frac {\Pi (x)}{\Pi (x^*)}\right\vert\geq$$
$$\geq\frac {n\delta_n^{2l}(x;x^*)}{C_{14}}\left\vert\frac {\Pi (x)}
{\Pi (x^*)}\right\vert.$$

Докажем оценку $(38)$. Пусть сначала $x^*-\pi\leq x<x^*-\frac
{\pi}{n}$, тогда
$$\delta_n(x;x^*)=\left\vert\frac {1}{n\sin\frac {x-x^*}{2}}\right\vert ;$$
в силу $(31)$, $(30)$ и $(36)$
\begin{gather*}
|\chi(x;x^*)-T(x)|=|T(x)|=\left\vert\int\limits_{x^*-\pi}^x
T'(u)du\right\vert\leq \\
\leq C_{13}n\int\limits_{x^*-\pi}^x\delta_n^{2(l-s)}(u;x^*)du=
C_{13}n\int\limits_{x^*-\pi}^x\left(\frac {1}{n\sin\frac {u-x^*}{2}}\right)
^{2(l-s)}du\leq \\
\leq\frac {C_{13}\pi^{2(l-s)}}{n^{2(l-s)-1}}\int\limits_{x^*-\pi}^x
\frac {du}{(u-x^*)^{2(l-s)}}\le \\
\le \frac {C_{13}\pi^{2(l-s)}}
{n^{2(l-s)-1}}\int\limits_{-\infty }^x\frac {du}{(u-x^*)^{2(l-s)}}= \\
=\frac {C_{13}\pi^{2(l-s)}}{(2(l-s)-1)n^{2(l-s)-1}}\left(\frac {1}
{(x^* -x)^{2(l-s)-1}}\right)\leq \\
\leq\frac {2C_{13}(\pi /2)
^{2(l-s)}}{(2(l-s)-1)n^{2(l-s)-1}}|\sin^{1-2(l-s)} \frac {(x^*-x)}2|=\\
=\frac {2C_{13}(\pi /2)^{2(l-s)}}{2(l-s)-1}\delta_n^{2(l-s)-1}
(x;x^*)=C_{19}\delta_n^{2(l-s)-1}(x;x^*). \end{gather*} Пусть
$x^*-\frac {\pi}{n}\leq x\leq x^*+\frac {\pi}{n}$. Тогда
$$\frac {2}{\pi}\leq\delta_n(x;x^*)\leq 1.$$
Поэтому с учётом $(36)$
\begin{gather*} |\chi (x;x^*)-T(x)|\leq
1+|T(x)| =1+\left\vert\int\limits_{x^*-\pi}^xT'(u)du
\right\vert\le \\
\leq 1+\left\vert\int\limits_{x^*-\pi}^{x^*-\frac {\pi}{n}}T'(u)du\right
\vert +\left\vert\int\limits_{x^*-\frac {\pi}{n}}^{x} T'(u)du\right\vert
\le \\
\le 1+C_{19}\delta_n^{2(l-s)-1}(x;x^*-\frac {\pi}{n})+\frac {2\pi C_{13}
n}{n}\leq 1+C_{19}+2\pi C_{13}=: \\
=:C_{15}\left(\frac {2}{\pi}\right)^{2(l-s)-1}\leq C_{15}\delta_n^{2(l-s)-1}
(x;x^*).
\end{gather*}

В случае $x^*+\frac {\pi}{n}\leq x\leq x^*+\pi$ рассуждаем так же, как и
в случае $x^*-\pi\leq x\leq x^*-\frac {\pi}{n}$, поскольку
$$\chi (x;x^*)-T(x)=1-T(x)=1-\int\limits_{x^*-\pi}^xT'(u)du=$$
$$=1-\int\limits_{x^*-\pi}^{x^*+\pi}T'(u)du+\int\limits_x^{x^*+\pi}T'(u)du=$$
$$=\int\limits_x^{x^*+\pi}T'(u)du.$$

В случае $x^*-2\pi\leq x\leq x^*-\pi$ имеем
\begin{gather*}
T(x)-\chi (x;x^*)=T(x)=\frac {1}{2\pi}x+\tilde {T}(x)=\\
=-1+\frac {1}{2\pi}(x+2\pi )+\tilde {T}(x+2\pi )=T(x+2\pi )-\chi
(x+2\pi ; x^*). \end{gather*}

Поскольку $x^*\leq x+2\pi\leq x^*+\pi$, то пользуясь предыдущими
двумя случаями, получаем оценку $(38)$.

Наконец, в случае $x^*+\pi\leq x\leq x^*+2\pi$ рассуждая
аналогично предыдущему случаю, получаем
\begin{gather*} |T(x)-\chi
(x;x^*)|=|T(x)-1|=\left\vert \frac 1{2\pi }x+\tilde T(x)-1
\right\vert =\\
=\left\vert \frac 1{2\pi }(x-2\pi )+\tilde T(x-2\pi )\right\vert
=|T(x-2\pi )|=\\
=|T(x-2\pi )-\chi (x-2\pi ;x^*)|\leq C_{15}\delta_n^{2(l-s)-1}(x;x^*).
\end{gather*}
Лемма доказана.

\medskip

\noindent $7^o$. Множество $O$.

Положим
$ O_i:=(x_{j+1}, x_{j-2})$, если $y_i\in [x_j, x_{j-1})$;
$$ O:=\bigcup\limits_{i\in\Bbb Z}O_i.                       \eqno (41)  $$
Будем писать $j\in H_*$, если $I_j\overline\in O$;
$j\in H$, если $j\in H_*$, $| j |\leq n$ и $j\ne -n$.
Обозначим $$\Bbb I:=[-\pi ,\pi ].$$

\medskip
\noindent $8^o$. ``Поправляющие'' функции.

В этом пункте введем ``поправляющие'' функции $T_j(x)$ и
$\overline {T}_j(x)$, первая из которых будет использована в
следующем пункте, а вторая -- в пункте $12^o$.

Положим
$$ x_j^*:=\frac {x_{j-1}+x_j}{2},             $$
$$c_{12}:=\max\{C_{12, s+2}, C_{12, s+3}\} ,  $$
$$n_1:=4(s+1)c_{12}n;$$
$$ x_j^{**}:=x_j^* + \frac {\pi}{n_1},        $$
$${\chi}_j(x):={\chi}_j(x,x_j).$$
Для каждого $j\in H$ обозначим
$$T_j(x):=\frac 12(T_{s+2,n_1}(x;x_j^*;s;Y)+T_{s+2,n_1}(x;x_j^{**};s;Y));
                                                               \eqno (42) $$
$$ \Pi_j(x):=\sin\left(\frac {x-x_j}{2}\right)\sin\left(\frac
{x-x_{j-1}}{2}\right)\Pi (x),                    $$
$$      Y_j:=Y\cup\{x_j\}\cup\{x_{j-1}\}; $$
$$ T_j^*(x):=T_{s+3,n_1}(x;x_j^*;s+1;Y_j)       $$
( т.е.
$$ T_j^*(x):=\frac{1}{d_{n_1}(x_j^*;s+1; Y_j)}\int\limits_{x_j^*-\pi}^x
J_{s+3,n_1}(x_j^* - t)\frac {\Pi_j(t)}{\Pi_j(x_j^*)}dt ); $$
$$ \overline {T}_j(x):=(T_j(x) - T_j^*(x)) sign\Pi (x_j).     \eqno (43) $$

\medskip

\noindent {\bf Л е м м а 6.} {\it Для каждого $j\in H$ функции
$T_j(x)$, $\overline {T}_j(x)$ имеют свойства
$$T_j(x)=\frac {1}{2\pi}x\ +\ \tilde {T}_j(x),           \eqno (44) $$
где $\tilde {T}_j\in{\Bbb T}_{(n_1-1)(s+2)+s}$;
$$\Pi (x_j^*)T_j'(x)\Pi(x)\geq 0,\ \ x\in\Bbb R,         \eqno (45) $$
$$|T_j'(x)|\geq c_{20}n_1\delta_{n_1}^{2(s+2)}(x;x_j^*)\left\vert\frac
{\Pi (x)}{\Pi (x_j^*)}\right\vert ,\ \ x\in\Bbb R,       \eqno
(46)
$$

$$|T_j'(x)|\geq c_{21}n_1\delta_{n_1}^{4(s+1)}(x;x_j^*),
\ \ x\in\Bbb R\setminus O,                               \eqno
(47)
$$
$$\sum\limits_{j\in H}| T_j(x) -\chi_j(x)|\leq c_{22},
\ \ x\in I,                                              \eqno
(48)
$$
$$\overline {T}_j\in{\Bbb T}_{(n_1-1)(s+3)+s+1},         \eqno (49) $$
$$|\overline {T}_j '(x)|\leq c_{23}n,\ \ x\in I_j,       \eqno (50) $$
$$\overline {T}_j '(x)\Pi (x)\geq 0,
\ \ x\in I\setminus I_j,                                 \eqno
(51)
$$
$$|\overline {T}_j '(x)|\geq |T_j '(x)|\geq c_{24}n\delta_{n_1}^{2(s+2)}
(x;x_j^*)\left\vert\frac {\Pi (x)}{\Pi (x_j^*)}\right\vert , \ \
x\in I_j\setminus I_j,                                 \eqno (52)
$$


$$|\overline T_j'(x)|\geq c_{25}n_1\delta_{n_1}^{4(s+1)}(x;x_j^*),
\ \ x\in I\setminus (O\cup I_j),                         \eqno
(53)
$$
$$\sum\limits_{j\in H}|\overline T_j(x)|\leq c_{26},\ \ x\in I.  \eqno (54) $$}

\noindent
{\it Д о к а з а т е л ь с т в о}. Поскольку $j\in H$, то
$$\min\limits_{i\in\Bbb Z}|x_j^*-y_i|>|x_j-x_j^*|=|x_{j-1}-x_j^*|=\frac
{\pi}{2n}=\frac {\pi 4(s+1)c_{12}}{2n_1}=$$
$$=2(s+1)c_{12}\frac {\pi}{n_1}>2sc_{12}\frac {\pi}{n_1},$$
$$\min\limits_{i\in\Bbb Z}|x_j^{**}-y_i|>2sc_{12}\frac {\pi}{n_1}.$$

Поэтому все соотношения леммы $6$, кроме $(47)$,$(48)$, $(53)$ и
$(54)$ следуют из леммы $5$. Неравенство $(47)$ немедленно
вытекает из $(46)$, поскольку для всех $y_i$ и $x\in\Bbb
R\setminus O$
$$\delta_{n}(x;x_j^*)\le 2\left\vert\frac {\sin\frac {x-y_i}2}{\sin
\frac {x_j^*-y_i}2}\right\vert.$$ По этой же причине неравенство
$(53)$ следует из $(52)$.

Так что докажем $(48)$, а затем $(54)$.

Представим разность
$$\alpha_j(x):=\chi (x;x_j^*)-T_j(x)$$ в виде
$$\alpha_j(x)=\frac {1}{2}(\chi (x;x_j^*)-T_{s+2,n_1}(x;x_j^*;s;Y)) +
\frac {1}{2}(\chi (x;x_j^{**})-                              $$
$$-T_{s+2,n_1}(x;x_j^{**};s;Y)) + \frac {1}{2}(\chi (x;x_j^*)
-\chi (x;x_j^{**})).                                         $$
Тогда в силу $(38)$ и оценки
$$\frac 14\delta_{n_1}(x;x_j^*)\leq\delta_{n_1}(x;x_j^{**})\leq 4\delta
_{n_1}(x;x_j^*),$$ имеем
\begin{gather*} |\alpha_j(x)|\leq\frac
{1}{2}c_{15}\delta_{n_1}^3(x;x_j^*)+\frac {64}2
c_{15}\delta_{n_1}^3(x;x_j^{**})+
\frac 12|\chi (x;x_j^*)-\chi (x;x_j^{**})|=\\
=\frac {65}2c_{15}\delta_{n_1}^3(x;x_j^*)
+\frac {1}{2}|\chi (x;x_j^*)-\chi (x;x_j^{**})|\leq \\
\leq\frac {65}2c_{15}\delta_{n_1}^3(x;x_j^*) +\frac
{1}{2}\delta_{n}^3(x;x_j^*)\leq c_{27}\delta_{n_1}^2(x;x_j^*),
\end{gather*}
 где $c_{15}=\max \{C_{15,s+2}, C_{15, s+3}\}$.

Поэтому
$$\sum\limits_{j\in H}|T_j(x)-\chi_j(x)|\leq c_{27}\sum\limits_{j=-n+1}^n
\delta_{n_1}^2(x;x_j^*)=:c_{27}\sum .$$
Оценим $\sum$. Пусть $x\in I_{\nu}$. Если $|x-x_j^*|<\pi$ и $j\ne\nu$, то
$$\delta_{n_1}(x;x_j^*)\leq\frac {1}{n_1|\sin (\frac {x-x_j^*}{2})|}\leq
\frac {\pi}{n_1|x-x_j^*|}<\frac 1{|j-\nu |};$$
если $\pi\le |x-x_j^*|<2\pi$, то воспользуемся  равенствами
$$\left\vert\sin\frac {x-x_j^*}2\right\vert=\left\vert\sin\frac {x-x_j^*}2
+\pi \right\vert=\left\vert\sin\frac {x-(x_j^*-2\pi)}2\right\vert=
\left\vert\sin\frac {x-x_{j+2s}^*}2\right\vert ,$$
и дело сводится к предыдущему случаю.
Поэтому учитывая равенство
$$\delta_{n_1}(x;x_j^*) = \delta_{n_1}(x;x_{j\pm 2n}^*)$$
получим
\begin{gather*} \sum
=\sum_{j=-n+1}^n\delta_{n_1}^2(x;x_j^*)
=\sum_{j=-n+1+\nu }^{n+\nu}\delta_{n_1}^2(x;x_j^*)=\\
=\delta_{n_1}^2(x;x_{\nu}^*) + \delta_{n_1}^2(x;x_{n+\nu}^*) +
\delta_{n_1}^2(x;x_{-n+1+\nu}^*)+\sum_{j=-n+2+\nu }^{n+\nu -1}\delta_{n_1}^2
(x;x_j^*)< \\
<3+\sum_{ j=-n+2+\nu, j\ne\nu}^{n+\nu -1}\frac 1{(j-\nu )^2}<3+2
\sum_{j=1}^{\infty}\frac 1{j^2}=3+\frac {\pi^2}3.
\end{gather*}
Так что $(48)$ доказано с $c_{22}=c_{27}(3 + \frac {\pi^2}3)$

Докажем $(54)$. Представим $\overline T_j(x)$ в виде
$$\overline T_j(x)sign\Pi (x_j)=-\alpha_j(x)+\chi (x;x_j^*)-T_j^*(x).$$
В силу $(38)$ имеем
$$|\overline T_j(x)|\le c_{27}\delta_{n_1}^2(x;x_j^*)+c_{15}\delta_{n_1}^2
(x;x_j^*)=(c_{27}+c_{15})\delta_{n_1}^2(x;x_j^*). $$
Поэтому
$$\sum\limits_{j\in H}|\overline T_j(x)|\leq (c_{27}+c_{15})\sum\limits
_{j=-n+1}^n\delta_{n_1}(x;x_j^*)\leq c_{26}.$$
Лемма доказана.

{\it З а м е ч а н и е  2.} Поскольку
$$\delta_{n_1}(x;x_j^*)>c_{27'}\frac 1{1+n\,dist\,(x,I_j)},$$
то $(53)$ влечёт
$$|\overline T_j'(x)|\ge c_{25}'n\left(\frac 1{1+n\,dist\,(x,I_j)}\right)^
{4(s+1)},\ \ \ \ x\in I\setminus (O\cup I_j) \eqno (53')$$

\medskip
\noindent $9^o$. Приближение функции $G_1(x)$.

\medskip
\noindent {\it О п р е д е л е н и е 3.} Через $N:=N(Y)>0$
обозначим число, имеющее следующее свойство : существуют две точки
$x^*$ и $x_*$ такие, что $\Pi (x)\geq 0$ при $x\in [x^* - \frac
{\pi}{N}, x^* + \frac {\pi}{N}]$, $\Pi (x)\leq 0$ при $x\in [x_* -
\frac {\pi}{N}, x_* + \frac {\pi}{N}]$.

\medskip
\noindent {\bf Л е м м а 7.} {\it Пусть непрерывно
дифференцируемая на $\Bbb R$ функция $G$ имеет вид
$$G (x) = Ax + \tilde {G}(x),$$
где $A=const$, $\tilde G$ - периодическая функция. Если
$$G'(x)\Pi (x)\geq 0,\ \ x\in\Bbb R,$$
и
$$\| G' \|\leq 1,$$
то при $n>N$ функция
$$V_n(x;G):=G(-\pi )+\sum\limits_{j\in H}(G(x_{j-1})-G(x_j))T_j(x)+$$
$$+\sum\limits_{j=-n+1, j\overline\in H}^n(G(x_{j-1}-G(x_j))(
\delta_j T_{s+2,n_1}(x;x^*;s;Y)+$$
$$+(1-\delta_j)T_{s+2,n_1}(x;x^*;s;Y)),$$
в которой
\[\delta_j:=\left\{\begin{array}{rl} 1,\ \mbox {если}\ \ \ \
G(x_{j-1})-G(x_j)\ge 0, \\
                   0,\ \mbox {если}\ \ \ \  G(x_{j-1})-G(x_j) <  0;
\end{array}\right.\]
имеет свойства
$$V_n(x;G)=Ax+\tilde V_n(x;G),                              \eqno (55) $$
где $\tilde V_n\in\Bbb T_{(s+2)(n_1-1)+s}$;
$$V_n'(x;G)\Pi (x)\geq 0,\ \ x\in\Bbb R;                    \eqno (56) $$
$$\|G-V_n\|\leq\frac {c_{28}\pi}n.                          \eqno (57) $$
} {\it Д о к а з а т е л ь с т в о.} Согласно $(44)$ и $(34)$
\begin{gather*} V_n(x;G)= G(-\pi )+\sum_{j\in
H}(G(x_{j-1})-G(x_j))\left(\frac x{2\pi}+
\tilde T_j(x)\right) + \\
+\sum\limits_{j=-n+1, j\overline\in H}^n
(G(x_{j-1}-G(x_j))\left(\frac x{2\pi} + \right .\\
\left .+\delta_j\tilde T_{s+2,n_1}(x;x^*;s;Y)+(1-\delta_j)\tilde
T_{s+2,n_1} (x;x_*;s;Y)\right)=:\\
=:\frac x{2\pi}\sum\limits_{j=-n+1}^n(G(x_{j-1})-G(x_j))+\tilde V_n(x;G)=\\
=\frac x{2\pi}(G(\pi )-G(-\pi ))+\tilde V_n(x;G)=Ax + \tilde V_n(x;G);
\end{gather*}
т.е. $(55)$ верно с учетом $(34')$.

Далее, если $j\in H$, то очевидно
$$sign (G(x_{j-1})-G(x_j))=sign\Pi(x_j^*),$$
поэтому с учетом $(45)$
$$sign\left((G(x_{j-1})-G(x_j))T_j'(x)\Pi (x)\right)=sign
(\Pi (x_j^*)T_j'(x)\Pi (x))\ge 0;$$ если $j\overline\in H$ и $sign
(G(x_{j-1})-G(x_j))>0$, то с учетом $(35)$
$$sign\left((G(x_{j-1})-G(x_j))\delta_jT_{s+2,n_1}'(x;x^*;s;Y)\Pi (x)\right)=
sign\Pi (x^*)\ge 0.$$ Если $j\overline\in H$ и  разность
$G(x_{j-1})-G(x_j)\le 0$, то
\begin{gather*}
sign\left((G(x_{j-1})-G(x_j))(1-\delta_j)T_{s+2,n_1}(x;x_*;s;Y)\Pi
(x)
\right)= \\
=-sign\Pi (x_*)\ge 0, \end{gather*}  т.е. неравенство $(56)$ тоже
верно. Осталось доказать оценку $(57)$. Для  этого рассмотрим
функцию
$$S_n(x;G)=G(-\pi )+\sum_{j=-n+1}^n(G(x_{j-1})-G(x_j))\chi_j(x).$$
Если $x\in (x_{\nu}, x_{\nu -1}]$, $\nu =\overline {-n+1, n}$, то
$$G(x)-S_n(x;G)=G(x)-G(x_{\nu})=\int_{x_{\nu}}^xG'(t)dt,$$
т.е.
$$\|G-S_n\|\leq\frac {\pi}n,$$
очевидно
$$|G(x_{j-1})-G(x_j)|\leq\frac {\pi}n,$$
Поэтому в силу $(48)$
$$\left\vert\sum_{j\in H}(G(x_{j-1})-G(x_j))(T_j(x)-\chi_j(x))\right\vert
\leq\frac {c_{22}\pi}n.$$
Теперь заметим, что количество индексов $j$ таких, что $j\overline\in H$
и $j=\overline {-n+1, n}$, не превосходит $6s$, тогда
$$\left\vert\sum\limits_{j=\overline {-n+1, n}, j\overline\in H}
(G(x_{j-1})-G(x_j))\chi_j(x)\right\vert\leq\frac {6s\pi}n,$$
аналогично с учётом $(38)$
\begin{gather*} \vert\sum\limits_{
j=\overline {-n+1, n}, j\overline\in H} (G(x_{j-1})-G(x_j))
(\delta_jT_{s+2, n_1}(x;x^*;s;Y)+\\
+(1-\delta_j)T_{s+2, n_1}(x;x_*;s;Y))\vert\le \\
\leq\left(|T_{s+2, n_1}(x;x^*;s;Y)|
+|T_{s+2, n_1}(x;x_*;s;Y))|\right)\frac {6s\pi}n\leq\\
\le (|T_{s+2, n_1}(x;x^*;s;Y)-\chi (x;x^*)|+|\chi (x;x^*)|+\\
+|T_{s+2, n_1}(x;x_*;s;Y)-\chi (x;x_*)|+|\chi (x;x_*)|)\frac {6s\pi}n
\leq\\
\leq (2+2c_{15})\frac {6s\pi}n.\end{gather*}

Таким образом находим
\begin{gather*}
\|G-V_n\|=\|(G-S_n)+\sum_{j\in H}(G(x_{j-1})-G(x_j))(\chi_j -T_j)+\\
+\sum\limits_{j=\overline {-n+1,n}, j\overline\in H}
(G(x_{j-1})-G(x_j))\chi_j -\\
-\sum\limits_{j=\overline {-n+1,n}, j\overline\in H}
(G(x_{j-1})-G(x_j))(\delta_j T_{s+2,n_1}(\ .\ ;x^*)+\\
+(1-\delta_j)T_{s+2,n_1}(\ .\ ;x_*))\|\leq\frac {\pi}n +\frac {c_{22}\pi}n+\\
+\frac {6s\pi}n + \frac {6s\pi}n(2+2c_{15}).\end{gather*}

Лемма доказана.

\medskip
Следствием леммы $7$ и соотношений $(24)$, $(17)$ является

\medskip
\noindent {\bf Л е м м а 8.} {\it При каждом $n>N$ существует
полином $\tilde V_n(x;G_1)$ такой, что
$$\tilde V_n\in\Bbb T_{(s+2)(n_1-1)+s}                       \eqno (55')$$
$$(B+\tilde V_n'(x;G_1))\Pi (x)\ge 0                         \eqno (56')$$
$$\|\tilde G_1-\tilde V_n\|\le\frac {c_2c_{28}\pi^r}{n^r}    \eqno (57')$$
} В самом деле, (55') и (56') очевидны, а (57') получается сразу,
как только мы заметим, что
$$\|G_1'\|\le c_2\left(\frac {\pi}n\right)^{r-1}.$$

\medskip
\noindent $10^o$. Определение функции $G_2(x)$.

\medskip
Если не существует ни одной пачки (опр. пачки см. п. $3^0$), то
$G_1(x)=f(x)$ и справедливость доказываемой теоремы вытекает из
леммы 8. Поэтому всюду далее предполагаем, что по крайней мере
одна пачка существует. В частности, $G_1(x)\not\equiv f(x)$.
Обозначим
$$G_2(x):=f(x)-G_1(x),                                    \eqno (58) $$
$$g_2(x):=G_2'(x).$$
Из $(16)$, $(19)$, $(24)$, $(15')$, $(8)$ и определения $W_2$
вытекают соответственно следующие свойства функции $G_2$:
$$G_2'\in c_{29}\Bbb W^{r-1}, \text {где}
\ \ c_{29}=c_4+1;                                        \eqno
(58')  $$
$$g_2(x)\Pi (x)\geq 0,\ \ x\in\Bbb R;                \eqno (58'') $$
$$G_2(x)=-Bx+\tilde G_2(x),                         \eqno (59)   $$
где $\tilde G_2\ -\ 2\pi$-периодическая функция;
$$G_2'(x)=0,\ \ x\in M^*;                                      \eqno (60)   $$
$$G_2'(x)=f'(x),\ \ x\in\Bbb R\setminus M^{**};                \eqno (61)   $$
$$|G_2'(x)|\geq h^{r-1}, x\in I_j,\ j\in V_2\setminus W_1.     \eqno (62)   $$

Поскольку функция $f$ периодическая, то без потери общности
будем считать $-n+1\in W_2\cap V_2$ и $n\in W_2\cap V_2$.

\medskip
\noindent $11^o$. Приближение без ограничений функции $G_2$.

\medskip
Обозначим
$$\Delta_t^r(G;x)=\sum_{j=0}^r(-1)^{r-j}\binom rj G(x+tj)$$
-- $r$-тую разность функции $G$ с шагом $t$ в точке $x$.

\medskip
\noindent {\bf Л е м м а 9.} {\it Пусть $G\in\Bbb W^r$, $r>1$, --
$2\pi$-периодическая функция ; $l>r+2$, $E\subset\Bbb R$, такое,
что если $y\in E$, то $(y\pm 2\pi )\in E$.  Если $$G''(x)=0,\ \
x\in\Bbb R\setminus E,$$ то тригонометрический полином
$$\Theta_{n,l,r}(x;G):=-(-1)^r\int\limits_{-\pi}^{\pi}
\sum\limits_{j=1}^{r}(-1)^{r-j}(_j^r)G(x+j\tau )J_{l,n}(\tau )d\tau $$
порядка
$$l(n-1)                                                  \eqno (63) $$
имеет свойства :

$$\| G-\Theta_{n,l,r} \|\leq c_{30}\left(\frac {\pi}n\right)^r;    \eqno (64) $$

$$| G'(x)-\Theta_{n,l,r}'(x,g) |\leq c_{31}\left(\frac {\pi}n\right)^{r-1}
\left(\frac {1}{1+ndist (x,E)}\right)^{2l-r},\ \ x\in\Bbb R, \eqno
(65) $$

$$\|G'-\Theta_{n,l,r}'\|\le c_{31}\left(\frac {\pi}n\right)^{r-1}; \eqno (66) $$

$$\|G''-\Theta_{n,l,r}''\|\le c_{32}\left(\frac {\pi}n\right)^{r-2}.\eqno (67) $$}

\noindent {\it Д о к а з а т е л ь с т в о}. Доказательство леммы
по сути ничем не отличается от рассуждений С.Б.Стечкина [22] (см.
также [23, с. 210]).  Для удобства читателя приведем эти
рассуждения.

Поскольку $G\in\Bbb W^r$, то (см., например  [7, с. 18])
$$|\bigtriangleup_t^r(x;G)|\leq |t|^r,\ \ x\in\Bbb R,          \eqno (68) $$
$$|\bigtriangleup_t^r(x;G')|=|\bigtriangleup_t^{r-1}(x;G')-
\bigtriangleup_t^{r-1}(x+t;G')|\leq 2|t|^{r-1},\ \ x\in\Bbb R;
                                       \eqno (69) $$
\begin{gather*}
|\bigtriangleup_t^r(x;G'')|=|\bigtriangleup_t^{r-2}(x;G'')-
2\bigtriangleup_t^{r-2}(x-t;G'') +\\
+\bigtriangleup_t^{r-2}(x-2t;G'')|\le 4|t|^{r-2},\ \ x\in\Bbb R
.\end{gather*}

1. Из условий леммы следует
$$G(x)-\Theta_{n,l,r}(x;G) =\int\limits_{-\pi}^{\pi}
\bigtriangleup_t^r(G;x)J_{l,n}(t)dt,\ \ x\in\Bbb R.             $$

В силу $(68)$ и $(29)$ получаем
$$\left\vert\int\limits_{-\pi}^{\pi}\bigtriangleup_t^r(G;x)J_{l,n}(t)dt
  \right\vert\le \int_{-\pi}^{\pi}|t|^rJ_{l,n}(t)dt\le
\frac {C_{12}\pi}{n^r},\ \ x\in\Bbb R.               $$ Тем самым
неравенство $(64)$ доказано.

\noindent
2. Аналогично пункту 1 имеем
$$\left\vert G'(x) - \Theta_{n,l,r}'(x;G)\right\vert =\left\vert\int
\limits_{-\pi}^{\pi}\bigtriangleup_t^r(G';x)J_{l,n}(t)dt\right\vert.
$$ И в силу $(69)$ и $(29)$ получаем
$$\left\vert\int\limits_{-\pi}^{\pi}\bigtriangleup_t^r(G';x)J_{l,n}(t)dt\right
\vert\leq\frac {2C_{12}}{n^{r-1}},\ \ x\in\Bbb R. $$ Тем самым
неравенство $(66)$ доказано а также доказано неравенство $(65)$ в
случае $x\in E$.

3. Докажем неравенство $(65)$ для $x\overline\in E$. Обозначим
$d=dist\ (x;E)$. Заметим, $\bigtriangleup_t^r(x;G')=0$ при
$|t|\leq\frac {d}{r}$ и разобьем разность
 $$\left\vert G'(x) - \Theta_{n,l,r}'(x;G)\right\vert =\left\vert\int
\limits_{-\pi}^{\pi}\bigtriangleup_t^r(G';x)J_{l,n}(t)dt\right\vert   $$
на два интеграла. Получим
$$\int\limits_{-\pi}^{\pi}\bigtriangleup_t^r(G';x)J_{l,n}
(t)dt=\int\limits_{-\pi}^{-d/r}\bigtriangleup_t^r(G';x)J_{l,n}(t)dt
+ \int\limits_{d/r}^{\pi}\bigtriangleup_t^r(G';x)J_{l,n}(t)dt. $$
Из $(28)$ и неравенства $(69)$ получаем оценку
\begin{gather*} |
G'(x) - \Theta_{n,l,r}'(x;G) |\leq\left\vert\int\limits_{-\pi}
^{-d/r}\bigtriangleup_t^r(G';x)J_{l,n}(t)dt\right\vert +\\
+ \left\vert\int
\limits_{d/r}^{\pi}\bigtriangleup_t^r(G';x)J_{l,n}(t)dt\right\vert
\leq\frac {4C_{11}}{n^{r-1}(1+nd/r)^{2l-r}}\le \\
\le C_{31}h^{r-1}\left(\frac
{1}{1+n\,dist(x,E)}\right)^{2l-r}.\end{gather*}

Тем самым доказано неравенство $(65)$ и в случае $x\overline\in
E$.

\noindent 4. Неравенство $(67)$ доказывается аналогично пункту 1,
или 2.

Лемма доказана.

\medskip
\noindent $12^o$. Приближение функции $G_2(x)$ на $\Bbb R\setminus
O$.

\medskip
\noindent
Обозначим
$$E:=\Bbb R\setminus M^*,\ \  C_{31}:=C_{31, s+1};$$
\begin{gather*} c_{34}:=\max \{ (4c_{29}c_{31})^{\frac 1{r+1}},\
\ \frac
{5^{4(s+1)}4c_{23}c_{29}c_{31}}{c_{25}'},  \\
\left(\frac {4c_{23}c_{29}c_{31}(1-(r-2)\pi )^{4(s+1)}}{c_{25}'}
\right)^{\frac 1{r-1}} \} , \hskip 3 cm (70') \end{gather*}
$$\ n_2:=c_{34}n.   $$

\medskip
\noindent {\bf Л е м м а 11.} {\it Функция
$$R(x):=\Theta_{n_2,2(s+1)+r,r}(x;\tilde G_2)-Bx+\frac {\pi^{r-1}}{2c_{23}n^r}
\sum\limits_{j\in V_2\setminus W_1}\overline T_j(x) \eqno (71) $$
имеет свойства :
$$R(x)=-Bx+\tilde R(x),                                    \eqno (72) $$
где $\tilde R\in\Bbb T_{(n_2-1)(2s+2 +r)}$;
$$\|\tilde R-\tilde G_2\|\leq\frac {c_{35}}{n^r};          \eqno (73) $$
$$R'(x)\Pi (x)\geq 0,\ \ x\in\Bbb R\setminus O.            \eqno (74) $$}

\noindent {\it Д о к а з а т е л ь с т в о}. Свойство $(72)$
следует из $(63)$ и $(49)$.  Неравенство $(73)$ следует из $(54)$,
$(57)$, $(64)$, и представления $(59)$. Докажем $(74)$. Обозначим
$$\frac {\pi^{r-1}}{2c_{23}n^r}\sum\limits_{j\in W_2\cap
V_2}\overline T_j(x) =:L(x) \eqno (75) $$ и заметим, что в силу
$(51)$
$$L'(x)\Pi (x)\geq 0                                        \eqno (76) $$
если $x\in I_{\nu}$, $\nu\overline\in V_2\cap W_2; $
$$L'(x)\Pi (x)\geq\frac {\pi^{r-1}}{2c_{23}n^r}\overline T_{\nu}'(x)\Pi (x),
                                                            \eqno (77) $$
если $\nu\in V_2\cap W_2$.

Без потери общности будем считать, что $x\in I=[-\pi ,\pi ]$.
Согласно $(53')$ имеем оценку
$$\frac {1}{n^r}|\overline T_j'(x)|\geq c_{25}'\frac {1}{n^{r-1}}\left(
\frac {1}{1+n dist\ (x,I_j)}\right)^{4(s+1)}, \eqno (78) $$
 $x\in I\setminus (O\cup I_j)$, $dist
(x,I_j)\leq\pi$. Обозначим
$$\Theta (x):=\Theta_{n_2, 2(s+1)+r,r}(x;\tilde G_2).$$
Из неравенства $(65)$ леммы $9$, $(58')$ и $(60)$ имеем
\begin{gather*} |\tilde G_2 '(x) - \Theta '(x)|\leq\frac
{c_{29}c_{31}\pi^{r-1}}{n_2^{r-1}}
\left(\frac {1}{1+n_2dist (x,E)}\right)^{4(s+1)}\leq\\
\le\frac {c_{29}c_{31}\pi^{r-1}}{n_2^{r-1}}\left(\frac 1{1+n dist
(x,E)} \right)^{4(s+1)}, \ \ x\in\Bbb R,  \end{gather*}  (79)
где, напомним, $E:=\Bbb R\setminus M^*$.

Из неравенств $(66)$ и $(58')$ имеем
$$\|\tilde G_2'-\Theta '\|\le\frac {c_{29}c_{31}\pi^{r-1}}{n_2^{r-1}}.
                                                                \eqno (80) $$
Производную $R'(x)$ представим в виде
$$R'(x)=\left(\Theta '(x)-\tilde G_2'(x)\right) +G_2'(x) +L'(x).\eqno (81) $$
Зафиксируем индекс $\nu\in H$, точку $x\in I_{\nu}$ и, считая для
определенности $\Pi (x)>0$, рассмотрим три случая. В случае 1. воспользуемся
оценкой
$$G_2'(x)\ge\left(\frac {\pi}n\right)^{r-1},$$
вытекающей из $(8)$, $(15')$ и $(58')$, а в случаях 2. и 3. --
оценкой $G_2'(x)\ge 0$.

1. Пусть $\nu\in V_2\setminus W_1$. В силу $(81)$, $(80)$, $(62)$,
$(58'')$, $(77)$, $(50)$ и $(70')$ имеем
$$R'(x)\ge -c_{29}c_{31}\left(\frac {\pi}{n_2}\right)^{r-1} +
\left(\frac {\pi}n\right)^{r-1} -\frac 12\left(\frac {\pi}n\right)^{r-1}=$$
$$=\left(\frac {\pi}n\right)^{r-1}\left(\frac 12 -\frac {c_{29}c_{31}}
{{c_{34}}^{r+1}}\right)>0.$$

2. Пусть теперь $\nu\in W_1$. Обозначим $j$ -- индекс ближайшего к $x$
 из отрезков $I_{\mu}$ таких, что $\mu\in V_2$ и $I_{\mu}\in\overline
{I\setminus M^{**}}$.

Тогда очевидно
$$dist\,(x, I_j)= dist\,(x, E)+\frac {\pi}n, $$
откуда
$$\frac 1{1+ndist\,(x, I_j)}>\frac 1{5(1+ ndist\,(x,E))}.    \eqno (82) $$
В силу $(81)$, $(79)$, $(76)$, $(78)$ и $(70')$ имеем
\begin{gather*}
R'(x)>\\
\frac {-c_{29}c_{31}}{n_2^{r-1}}\left(\frac 1{1+ndist\,(x,E)}\right)
^{4(s+1)}+
\frac {c_{25}'\pi^{r-1}}{2c_{23}n^{r-1}}\left(\frac 1{1+ndist\,(x,I_j)}
\right)^{4(s+1)}>\\
>\left(\frac {\pi}n \right)^{r-1}\left( \frac {-c_{29}c_{31}}{c_{34}} +
\left( \frac 15 \right)^{4(s+1)}\frac {c_{25}'}{2c_{23}}\right)
\left( \frac 1{1+ndist\,(x, E)} \right)^{4(s+1)}>0.
\end{gather*}

3. Наконец, пусть $\nu\in W_2\setminus V_2$. Обозначим $j$--
индекс ближайшего к $x$ из отрезков $I_{\mu}$, таких, что $\mu\in
V_2$. В силу Леммы 2
$$dist (x, I_j)\le\frac {\pi (r-2)}n.                         \eqno (83) $$
В силу $(81)$, $(80)$, $(76)$, $(78)$, $(83)$ и $(70')$ имеем
\begin{gather*} R'(x)>-c_{29}c_{31}\left(\frac
{\pi}{n_2}\right)^{r-1}+\frac {c_{25}'}{2c_{23}} \left(\frac
{\pi}n\right)^{r-1}\left(\frac 1{1+n\,dist(x,I_j)}\right)^{4(s+1)}
>\\
>\left(\frac {\pi}n\right)^{r-1}\left(\frac {-c_{29}c_{31}}{c_{34}^{r-1}}+
\frac {c_{25}'}{2c_{23}(1+(r-2)\pi )^{4(s+1)}}\right)>0.
\end{gather*} Лемма доказана.

$13^0$. Приближение функции на О.

Обозначим через $N_1:=N_1(Y)$ число такое, что при всех $n>N_1$ и
любом $i$
$$y_{i-1}-y_i>\frac {6\pi}n.                                    \eqno (84) $$
Зафиксируем точку $y_i\in\Bbb I$. Пусть $y_i\in [x_j, x_{j-1})$.
Положим
$$\overline Y_i:=Y\setminus\{y_i\}\cup\{x_{j+1}\},\ \
  \underline Y_i:=Y\setminus\{y_i\}\cup\{x_{j-2}\};             \eqno (85) $$
$$\overline\Pi_i(x):=\sin\frac {(x-x_{j+1})}2\prod\limits_{\nu =0, \nu\ne
i}^{2s-1} \sin\frac {(x-y_i)}2, \eqno (86) $$
$$\underline\Pi_i(x):=\sin\frac {(x-x_{j-2})}2\prod\limits_{\nu
=0, \nu\ne i}^{2s-1} \sin\frac {(x-y_i)}2; \eqno (87) $$
$$\overline P_i(x):=\frac 12(T_{s+2,n_1}(x;x_{j-2}^*,s,\overline Y_i) +
T_{s+2,n_1}(x;x_{j-2}^{**},s,\overline Y_i)), \eqno (88) $$
$$\underline P_i(x):=\frac 12(T_{s+2,n_1}(x;x_{j+2}^*,s,\underline Y_i) +
T_{s+2,n_1}(x;x_{j+2}^{**},s,\underline Y_i)); \eqno (89) $$
$$\overline K_i(x):=(T_{j+2}(x)-\overline P_i(x))sign\Pi (x_{j+2}^*),
                                                                 \eqno (90) $$
$$\underline K_i(x):=(T_{j-2}(x)-\underline P_i(x))sign\Pi (x_{j-2}^*);
                                                                 \eqno (91) $$
$$U_i(x):=\frac 1{n^r}(T_{j+2}(x)-T_{j-2}(x))sign\Pi (x_{j+2}^*).\eqno (92) $$
Из $(44)$ и $(34)$ следует
$$\overline K_i, \underline K_i, U_i \in\Bbb T_{(s+2)(n_1-1)+s}.   \eqno (93) $$
Из $(45)$ и $(35)$ для $x\in\Bbb I\setminus O_i$ следует
\begin{gather*}
\overline K_i'(x)sign\,\Pi (x)=T_{j+2}'(x)sign\Pi (x)sign\Pi (x_{j+2})-\\
-\overline P_i'(x)sign\Pi (x)sign\Pi (x_{j+2})=\\
=T_{j+2}'(x)sign\Pi (x)sign\Pi (x_{j+2})+ \overline
P_i'(x)sign\overline\Pi_i (x)sign\overline\Pi (x_{j-2})\ge 0, \ \
\ \ \ \ \ \ \ \ \ \ \ (94)
\end{gather*}
\begin{gather*}
\underline K_i'(x)sign\,\Pi (x)=T_{j-2}'sign\Pi (x)sign\Pi (x_{j-2})- \\
-\underline P_i'(x)sign\Pi (x)sign\Pi (x_{j-2})=\\
=T_{j-2}'(x)sign\Pi (x)sign\Pi (x_{j-2})+ \underline
P_i'(x)sign\underline\Pi_i (x)sign\underline\Pi (x_{j+2})\ge 0,
\hskip 3 cm (95)         \end{gather*}
\begin{gather*} \overline K_i'(y_i)=-\overline P_i'(y_i)sign\Pi
(x_{j+2}^*)=
\overline P_i'(y_i)sign\overline\Pi_i (x_{j-2})=\\
=\left(\overline P_i'(y_i)sign\overline\Pi_i(x_{j-2})sign\overline\Pi_i
(y_i)\right)sign\overline\Pi_i(y_i)=\\
=\left\vert\overline P_i'(y_i)\right\vert sign\overline\Pi_i(y_i)
=\left\vert\overline P_i'(y_i)\right\vert sign\Pi(x_{j-2}).\hskip
7 cm (96)
\end{gather*}
\begin{gather*} \underline K_i'(y_i)=-\underline P_i'(y_i)sign\Pi
(x_{j-2}^*)=
\underline P_i'(y_i)sign\underline\Pi_i (x_{j+2})=\\
=\left(\underline P_i'(y_i)sign\underline\Pi_i (x_{j+2})sign\underline
\Pi_i (y_i)\right) sign\underline\Pi_i (y_i)=\\
=\left\vert\underline P_i'(y_i)\right\vert sign\underline\Pi_i
(y_i)= \left\vert\underline P_i'(y_i)\right\vert sign\Pi
(x_{j+2}). \hskip 7 cm  (97)
\end{gather*}
В силу $(37)$
\begin{gather*} \left\vert \overline P_i'(y_i)
\right\vert\ge\frac n{2C_{14,s+2}}\Bigg\lgroup
\delta_n^{2(s+2)}(y_i;x_{j-2}^*)\left\vert\frac
{\overline\Pi_i(y_i)}
{\overline\Pi_i(x_{j-2}^*)}\right\vert+\\
+\delta_n^{2(s+2)}(y_i;x_{j-2}^{**})\left\vert \frac
{\overline\Pi_i(y_i)} {\overline\Pi_i(x_{j-2}^{**})}\right\vert
\Bigg\rgroup \ge c_{37}n,  \hskip 8 cm (98)
\end{gather*}
аналогично $$\left\vert\underline P_i'(y_i)\right\vert\ge c_{37}n.
\eqno (99)$$

Из $(38)$ для всех $x\in\Bbb I$ следует
\begin{gather*}
\left\vert\overline K_i(x)\right\vert =\left\vert T_{j+2}(x)-
\chi(x;x_{j+2}^*)-\frac 12(T_{s+2,n_1}(x;x_{j-2}^*;s;Y)-\right .\\
\chi(x;x_{j-2}^*))
-\frac 12(T_{s+2,n_1}(x;x_{j-2}^{**};s;Y)-\chi(x;x_{j-2}^{**}))+\\
\left . +\left(\chi(x;x_{j+2}^*)-\frac 12\chi(x;x_{j-2}^*)
-\frac 12\chi(x;x_{j-2}^{**})\right)\right\vert \le \\
\le 2C_{15,s+2}+2=:C_{38}. \hskip 8 cm (100)  \end{gather*}

Аналогично, $$\left\vert\underline K_i(x)\right\vert\le C_{38}.
\eqno (101)$$

Наконец положим \[ K_i(x):=\left\{\begin{array}{rl}
\left\vert\frac {R'(y_i)}{\overline K_i'(y_i)}\right\vert
\overline K_i(x), &\
\mbox {если}\ \ R'(y_i)sign\Pi (x_{j+2})\ge 0,\\
          \left\vert\frac {R'(y_i)}{\underline K_i'(y_i)}\right\vert
\underline K_i(x),&\ \mbox {если}\ \ R'(y_i)sign\Pi (x_{j+2})< 0.
             \end{array} \hskip 5 cm  (102) \right.\]

{\bf Л е м м а 12.} {\it Имеют место соотношения
$$|R'(x)+K_i'(x)-G_2'(x)|\le c_{39}\frac {|x-y_i|}{n^{r-2}},\ \ x\in O_i,
                              \eqno (103) $$
$$\|K_i\|\le\frac {c_{40}}{n^r}                             \eqno (104) $$
$$K_i'(x)\Pi (x)\ge 0,\ \ x\in\Bbb I\setminus O_i           \eqno (105) $$}

{\it Д о к а з а т е л ь с т в о.} Сначала докажем равенство
$$ R'(y_i)+K_i'(y_i)-G_2'(y_i)=0.                           \eqno (106) $$
Действительно, по условию $G_2'(y_i)=0$. Если
$$R'(y_i)sign\Pi (x_{j+2})\ge 0,$$
то в силу $(102)$ и $(96)$
\begin{gather*}
R'(y_i)+K_i'(y_i)=R'(y_i)+\left\vert\frac {R'(y_i)}
{\overline K_i'(y_i)}\right\vert\overline K_i'(y_i)=\\
=R'(y_i)+|R'(y_i)|sign\Pi (x_{j-2})=R'(y_i)-|R'(y_i)|sign\Pi (x_{j+2})=0.
\end{gather*}
Если $R'(y_i)sign\Pi (x_{j+2})<0$, то в силу $(102)$ и $(97)$

\begin{gather*} R'(y_i)+K_i'(y_i)=R'(y_i)+\left\vert\frac
{R'(y_i)}{\underline K_i'(y_i)}
\right\vert\underline K_i'(y_i)=\\
=R'(y_i)+\left\vert R'(y_i)\right\vert sign\Pi (x_{j+2})=0.
\end{gather*}
Таким образом равенство $(106)$ доказано. Теперь заметим, что в
силу $(73)$ и неравенства Бернштейна
\begin{gather*}
|R'(y_i)|=|R'(y_i)-G_2'(y_i)|=|\tilde R'(y_i)-\tilde G_2'(y_i)|\le\\
\le \|\tilde R'-\tilde G_2'\|\le c_{41}n\|\tilde R-\tilde
G_2\|\le\frac {c_{42}}{n^{r-1}}. \end{gather*}
С другой стороны,
согласно $(96)$ -- $(99)$
$$\left\vert\overline K_i'(y_i)\right\vert >c_{37}n,\ \ \
\left\vert\underline K_i'(y_i)\right\vert >c_{37}n. $$ Поэтому, с
учётом $(100)$ и $(101)$
$$\|K_i\|\le\frac {C_{38}c_{42}}{c_{37}n^r}=:\frac {c_{40}}{n^r}.$$
Таким образом доказано неравенство $(104)$. Неравенство $(73)$,

\noindent $(104)$ и неравенство Бернштейна немедленно влекут
оценку
$$\|R''+K_i''-G_2''\|\le\frac {c_{43}}{n^{r-2}}, $$
которая вместе с $(106)$ даёт
\begin{gather*}
|R'(x)+K_i'(x)-G_2'(x)|=|R'(x)+K_i'(x)-G_2'(x)-R'(y_i)-K_i'(y_i)+G_2'(y_i)|
\\
=\left\vert\int\limits_{y_i}^x(R''(t)+K_i''(t)-G_2''(t))dt\right\vert\le
\left\vert\int\limits_{y_i}^x\frac
{c_{43}}{n^{r-2}}dt\right\vert=\frac
{c_{39}|x-y_i|}{n^{r-2}}.\end{gather*}

Таким образом и неравенство $(103)$ леммы $12$ доказано.
Соотношение $(105)$ есть очевидное следствие $(102)$, $(94)$ и
$(95)$. Лемма доказана.

Напомним (см. $(92)$)
$$U_i(x):=\frac 1{n^r}(T_{j+2}(x)-T_{j-2}(x))sign\Pi (x_{j+2}^*).$$
Согласно $(45)$
$$U_i'(x)\Pi (x)\ge 0,\ \ x\in\Bbb I.$$
Поскольку при $x\in O_i$
$$\delta(x;x_{j+2}^*)>c_{44}, \ \ \delta(x;x_{j-2}^*)>c_{44},$$
$$\left\vert\frac {\Pi (x)}{\Pi (x_{j+2}^*)}\right\vert >c_{45}|x-y_i|n,$$
и
$$\left\vert\frac {\Pi (x)}{\Pi (x_{j-2}^*)}\right\vert >c_{45}|x-y_i|n,$$
то в силу $(46)$
$$\left\vert U_i'(x)\right\vert\ge\frac {c_{46}|x-y_i|}{n^{r-2}},
\ \ x\in O_i.                                               \eqno
(107) $$ Кроме того, согласно $(38)$
$$\| U_i\| =\frac 1{n^r}\left\Vert\left(T_{j+2}-\chi_{j+2}
\right) +\left(\chi_{j-2}-T_{j-2}\right)
+\chi_{j+2}-\chi_{j-2}\right\Vert \le\frac {c_{47}}{n^r}$$ и
согласно $(44)$
$$U_i\in\Bbb T_{(n_1-1)(s+2)+s}.                              \eqno (108) $$
Таким образом, справедлива

\medskip
{\bf Л е м м а 13.} {\it Имеют место соотношения
$$\| U_i \|\le\frac {c_{47}}{n^r},                            \eqno (109) $$
$$\left(R'(x)+K_i'(x)+\frac {c_{39}}{c_{46}}U_i'(x)\right)\Pi(x)\ge 0,
                 x\in\Bbb I.                                  \eqno (110) $$
$$U_i\in\Bbb T_{(n_1-1)(s+2)+s},                              \eqno (111) $$
$$U_i'(x)\Pi (x)\ge 0,\ \ x\in\Bbb I,                         \eqno (112) $$ }
{\it Д о к а з а т е л ь с т в о.} Соотношения $(111)$, $(109)$ и
$(112)$ суть соотношения $(108)$, $(104)$ и $(105)$. Неравенство
$(110)$ следует из $(112)$, $(107)$ и $(103)$. Действительно,
\begin{gather*}
\left(R'(x)+K_i'(x)+\frac {c_{39}}{c_{46}}U_i'(x)\right)\Pi (x)=\\
=\left(R'(x)+K_i'(x)+\frac {c_{39}}{c_{46}}U_i'(x)-G_2'(x)\right)\Pi (x)+
G_2'(x)\Pi (x)\ge \\
\ge \frac {c_{39}}{c_{46}}U_i'(x)\Pi (x)+(R'(x)+K_i'(x)-G_2'(x))\Pi (x)\ge
\\
\ge \frac {c_{39}}{c_{46}}|U_i'(x)||\Pi (x)|+(R'(x)+K_i'(x)-G_2'(x))\Pi (x)
\ge \\
\ge\frac {c_{39}|x-y_i|}{n^{r-2}}|\Pi(x)|-\frac
{c_{39}|x-y_i|}{n^{r-2}} |\Pi(x)|=0. \end{gather*}
   Лемма доказана.

Наконец, из предыдущих лемм вытекает

\medskip
{\bf Л е м м а 14.} {\it При $n>N_2:=\max\{ N, N_1\}$ полином
$$\tau (x)=V_n(x;G_1)+R(x)+\sum\limits_{i=0}^{2s-1}\left(K_i(x)+\frac
{c_{39}}{c_{46}}U_i(x)\right)$$
порядка $(n-1)(s+2)+s$ имеет свойства
$$\tau'(x)\Pi (x)\ge 0,\ \ x\in\Bbb R,$$
$$\|f-\tau \|\le\frac c{n^r}.$$}
{\it Д о к а з а т е л ь с т в о.} В самом деле, первое
неравенство леммы следует из соотношений (56) леммы 7, (74) леммы
11 и (110) леммы 13. Второе неравенство леммы следует из
соотношений (57) леммы 7, (73) леммы 11, (104) леммы 12 и (109)
леммы 13. Лемма доказана.

\medskip
Таким образом при $n>N_2(Y)$ теорема доказана. Для $n\le N_2(Y)$ теорема
является следствием частного случая $n=0$, который вытекает из неравенства
Уитни.

\bigskip

\begin{center}
ЛИТЕРАТУРА \end{center}

\noindent 1. Pleshakov M.G. {\it Comonotone Jackson's Inequality,}
J. Approx. Theory,  {\bf 99} (1999), 409--421.

\noindent 2. Гилевич Я., Шевчук И.А. {\it Комонотонное
приближение}, Фунд. и прикл. математика, {\bf 2}

(1996), 2, 319--363.

\noindent 3. Дзюбенко Г.А. {\it Поточечная оценка комонотонного
приближения,} УМЖ,  {\bf 46} (1994), 11,

1467--1472.

\noindent 4. Никольский С.М. {\it О наилучшем приближении
многочленами  функций, удовлетворя-

ющих условию Липшица,} Изв. АН СССР, сер. мат., {\bf 10} (1946),
4, 295--317.

\noindent 5. Шевчук И.А. {\it О коприближении монотонных функций,}
ДАН СССР, {\bf 308} (1989), 3,

 537--541.

\noindent 6. Шевчук И.А. {\it Приближение монотонных функций
монотонными многочленами,} Ма-

тем. сборник., {\bf 183} (1992), 5, 63--78.

\noindent 7. Шведов А.С. {\it Теорема Джексона в $L_p$, $0<p<1$,
для ал\-ге\-б\-ра\-и\-че\-ских многочленов и

порядки комонотонных приближений}, Мат. зам., {\bf 25} (1979), 1,
107--117.

\noindent 8. Шевчук И.А. {\it Приближение многочленами и следы
непрерывных на отрезке функций. }

``Наукова думка'', Киев, 1992, 225 стр.

\noindent 9. Шведов А.С. {\it Комонотонное приближение функций
многочленами,} ДАН СССР, {\bf 250}

(1980), 1, 39--42.

\noindent 10. Шведов А.С. {\it Коприближение кусочно-монотонных
функций многочленами,} Мат. за-

метки,  {\bf 30} (1981), 6, 839--846.

\noindent 11. Beatson R.K., Leviatan D. {\it On comonotone
approximation,} Canad. Math.  Bull. {\bf 26} (1983),

220--224.

\noindent 12. DeVore R.A. {\it Degree of Approximation,}
Approximation theory II, (Lorentz G.G., Chui C.K.,

Schumaker L. eds.), New York: Acad. Press, 1976, 117--162.

\noindent 13. DeVore R.A. {\it Monotone Approximation by
Polynomials,} SIAM J. Math. Anal., {\bf 8} (1977),  5,

906--921.

\noindent 14. DeVore R.A., Yu X.M. {\it Pointwise Estimates for
Monotone Polynomial Approximation,}

Constr. Approx.,  {\bf 1} (1985), 4, 323--331.

\noindent 15. Leviatan D. {\it Pointwise Estimates for Convex
Polynomial Approximation.} Proc. Amer.

Math. Soc.,  {\bf 98} (1986), 471--474.

\noindent 16. Leviatan D., Shevchuk I.A., {\it Some Positive
Results and Counterexamples in Comonotone

Approximation}, J. Approx. Theory, {\bf 89} (1997), 195-206.

\noindent 17. Newman D.J., Passow E., Raymon L. {\it Piecewise
monotone polynomial approximation,}

Trans. Am. Math. Soc., {\bf 172} (1972), 465--472.

\noindent 18. Passow E., Raymon L., Roulier J.A. {\it Comonotone
polynomial approximation,} J.Approx.Theory,

{\bf 11} (1974), 221--224.

\noindent 19. Passow E., Raymon L. {\it Monotone and comonotone
approximation,} Proc. Amer. Math. Sos.,

{\bf 42} (1974), 390--394.

\noindent 20. Shevchuk I.A. {\it Whitney's Inequality and
Coapproximation}, East Journal on Approximations,

{\bf 1} (1995), 4 , 479--500.

\noindent 21. Плешаков М.Г. {\it О кусочно-монотонном приближении
периодических функций три-

гонометрическими полиномами}, Теор. основы и конструирование
численных алгорит-

мов решения задач мат. физики, тез. докл. 11-ой Всеросс. конфер.,
Пущино, 1996.

\noindent 22. Стечкин С.Б. {\it О наилучшем приближении
периодических функций тригонометриче-

скими полиномами}, Доклады АН СССР,  {\bf 83} (1952), 5, 651--654.

\noindent 23. Дзядык В.К. {\it Введение в теорию равномерного
приближения функций полиномами.}

''Наука'', M., 1977, 512 стр.

\end{document}